\documentclass[12pt]{article}
\usepackage[cp1251]{inputenc}
\usepackage[russian]{babel}
\usepackage{amssymb}
\usepackage{amsmath}

\begin{document}

\begin{center}
\textbf{\large On $L^2$ -functions with bounded spectrum}
\end{center}

\begin{center}
Vladimir Lebedev
\end{center}

\begin{quotation}
{\small \textsc{Abstract.} We consider the class $PW(\mathbb R^n)$
of functions in $L^2(\mathbb R^n)$, whose Fourier transform has
bounded support. We obtain a description of continuous maps
$\varphi : \mathbb R^m\rightarrow\mathbb R^n$ such that
$f\circ\varphi\in PW(\mathbb R^m)$ for every function $f\in
PW(\mathbb R^n)$. Only injective affine maps $\varphi$ have this
property.

  Bibliography: 5 items.

  Keywords: Fourier transform, spectrum bounded functions,
superposition operators.

  AMS 2010 Mathematics Subject Classification. Primary 42B10.
Secondary 30D15}
\end{quotation}

\quad

  Let $PW(\mathbb R^n)$ be the class of functions $f$ of the form
$$
f(t)=\int_B g(u)e^{i(u, t)} du, \qquad t\in\mathbb R^n,
$$
where $g$ is an arbitrary function in $L^2(\mathbb R^n)$ and
$B\subset\mathbb R^n$ is an arbitrary ball ($(u, t)$ stands for
the inner product). For $n=1$ the corresponding class was
considered by Paley and Wiener [1]. We note that functions of
class $PW$ appear in the problems related to signal processing and
are often called bandlimited signals (in this connection see,
e.g., the references in [2]).

  It is clear that a function $f$ is in $PW(\mathbb R^n)$
if and only if it is continuous, belongs to $L^2(\mathbb R^n)$,
and has bounded spectrum, that is, its Fourier transform
$\widehat{f}$ vanishes outside a certain ball. For each such
function $f$ we have $\widehat{f}\in L^1(\mathbb R^n)$ and
$$
f(t)=\int_{\mathbb R^n}\widehat{f}(u)e^{i(u, t)}du,
\qquad t\in\mathbb R^n.
$$

  Let $\varphi : \mathbb R^n\rightarrow\mathbb R^n$ be a
non-degenerate affine map. It is easy to see that then for every
function $f\in PW(\mathbb R^n)$ the superposition
$(f\circ\varphi)(t)=f(\varphi(t))$ is also in $PW(\mathbb R^n)$.
It suffices to note that if $\varphi(t)=At+b$, where $A$ is an
invertible $n\times n$ matrix and $b\in\mathbb R^n$, then
$$
|\widehat{f\circ\varphi}(u)|=|\det A|^{-1}|\widehat{f}((A^{-1})^*
u)|
\eqno(1)
$$
for every function $f\in L^2(\mathbb R^n)$. Here $A^{-1}$ is the
inverse of the matrix $A$ and $(A^{-1})^*$ is the transpose of the
matrix $A^{-1}$.

  In [2] Azizi, Cochran, and McDonald showed that if
$\varphi$ is a homeomorphism of the line $\mathbb R$ onto itself
such that for every $f\in PW(\mathbb R)$ we have $f\circ\varphi\in
PW(\mathbb R)$, then the map $\varphi$ is affine. The same authors
posed the question [3] as to whether the similar assertion holds
in the multidimensional case. In the present work we shall show,
in particular, that the answer to this question is positive,
moreover this is true not only for homeomorphisms but for
arbitrary continuous maps $\varphi : \mathbb R^n\rightarrow\mathbb
R^n$. Actually we shall give the description of continuous maps
$\varphi : \mathbb R^m\rightarrow\mathbb R^n$ which have the
property that $f\circ \varphi\in PW(\mathbb R^m)$ whenever $f\in
PW(\mathbb R^n)$.

  Let $\varphi : \mathbb R^m\rightarrow\mathbb R^n$ be an affine
map. We have $\varphi(t)=At+b$, where $A$ is an $n\times m$ matrix
and $b\in\mathbb R^n$. It is clear that the map $\varphi$ is
injective if and only if the kernel
$$
\mathrm{ker} A=\{x\in\mathbb R^m : Ax=0\}
$$
of the matrix $A$ is trivial, that is, $\mathrm{ker} A=\{0\}$.

\quad

   \textbf{Theorem.} \emph{Let $\varphi$ be a continuous map
of $\mathbb{R}^m$ into $\mathbb{R}^n$. The following conditions
are equivalent: \emph{(i)} for every function $f\in
PW(\mathbb{R}^n)$ the superposition $f\circ\varphi$ belongs to
$PW(\mathbb{R}^m)$; \emph{(ii)} $\varphi$ is an injective affine
map.}

\quad

   Note the immediate consequence of this theorem: for
$n<m$ there are no continuous maps $\varphi : \mathbb
R^m\rightarrow\mathbb R^n$ such that for every $f\in
PW(\mathbb{R}^n)$ the superposition $f\circ\varphi$ is in
$PW(\mathbb{R}^m)$.

\quad

   \emph{Proof of the theorem.} We shall first prove the
implication (i)$\Rightarrow$(ii). Recall that an entire (i.e. an
analytic in the complex plane $\mathbb C$) function $f$ is said to
be of exponential type if $|f(z)|=O(e^{b|z|}), ~z\in\mathbb C,$
where $b$ is a positive constant. It is well known that for each
such function $f$ with zeros $\{z_k, ~k=1, 2, \ldots\}$ we have
$$
f(z)=e^{c_0+c_1z} \Pi(z),
$$
where
$$
\Pi(z)=z^m \prod_{k :\, z_k\neq 0} \bigg(1-\frac{z}{z_k}\bigg)e^{z/z_k}
$$
is the canonical product constructed with the zeros $\{z_k\}$ of
$f$ counting their multiplicities (this is a particular case of
the Hadamard factorization theorem, see, e.g. [4, Ch. VIII]).

\quad

  \textbf{Lemma 1.} \emph{Let $A(z), ~B(z), ~C(z), ~D(z)$ be
entire functions. Assume that $C$ and $D$ are of exponential type.
Assume that the functions $A$ and $C$ are not identically zero.
Suppose that $\psi$ is a continuous real function on $\mathbb{R}$
such that for all real $x$ we have
$$
A(x)\psi(x)=B(x), \eqno (2)
$$
$$
C(x)e^{i\psi(x)}=D(x). \eqno (3)
$$
Then $\psi$ is linear \emph{(}affine\emph{)} function,
$\psi(x)=ax+b, ~a, b\in \mathbb{R}$.}

\quad

   \emph{Proof.} Obviously (see (3)) the function
$D(z)$ is not identically zero. Denote by $\mathcal{N}$ the union
of the sets of zeros of the functions $A(z), C(z), D(z)$. All
points of the set $\mathcal{N}\subset\mathbb C$ are isolated.

   From (2) and (3) we obtain that
$$
e^{iB(x)/A(x)}=\frac{D(x)}{C(x)}
$$
for all $x\in \mathbb R \setminus\mathcal{N}$. Hence (by
uniqueness theorem, see, e.g. [5, Ch. III, \S ~6]) for all $z\in
\mathbb C\setminus\mathcal{N}$ we have
$$
e^{iB(z)/A(z)}=\frac{D(z)}{C(z)}, \eqno (4)
$$
$$
e^{-iB(z)/A(z)}=\frac{C(z)}{D(z)}. \eqno (5)
$$

   All singular points of the functions $B/A, ~D/C, ~C/D$
are either removable or poles. From (4) it follows that each pole
of $B/A$ if there any should be an essential singular point of
$D/C$, which is impossible. So $B/A$ has only removable singular
points. Then by (4), (5) the singular points of the functions
$D/C$ and $C/D$ are removable as well and hence, the zeros of $C$
and $D$ coincide counted with their multiplicities. Denote by
$\Pi$ the canonical product constructed with these zeros. Taking
into account that $C$ and $D$ are of exponential type we obtain
$C(z)=e^{\gamma(z)}\Pi(z)$ and $D(z)=e^{\delta(z)}\Pi(z),
~z\in\mathbb C,$ where $\gamma$ and $\delta$ are linear functions
of $z$. Relation (4) yields
$$
e^{iB(z)/A(z)}=\frac{D(z)}{C(z)}=e^{\delta(z)-\gamma(z)}
$$
for all $z\in \mathbb C\setminus\mathcal{N}$. So, for $x\in
\mathbb R\setminus\mathcal{N}$ we have (see (2))
$$
e^{i\psi(x)}=e^{il(x)}, \eqno (6)
$$
where $l(x)=-i(\delta(x)-\gamma(x))$ is a linear function. Since
the function $\psi$ is continuous and the set $\mathcal{N}$ is at
most countable, relation (6) holds for all $x\in\mathbb R$.
Therefore, $\psi(x)=l(x)+2\pi k(x), ~x\in\mathbb R$, where $k(x)$
takes only integer real values. Due to continuity of $\psi$ the
function $k(x)$ is constant. The lemma is proved.

\quad

  The following lemma is trivial.

\quad

  \textbf{Lemma 2.} \emph{Let $f\in PW(\mathbb{R}^n)$.
Let $l : \mathbb R\rightarrow\mathbb R^n$ be an affine map. Then
there exists an entire function $F(z)$ of exponential type such
that $F(x)=f\circ l(x)$ for all $x\in\mathbb{R}$.}

\quad

   \emph{Proof.} We have $l(x)=a+xb, ~x\in\mathbb R$,
where $a$ and $b$ are certain vectors in $\mathbb R^n$. Let
$B=B(0, r)$ be a ball in $\mathbb R^n$ centered at $0$ and of
radius $r$ that contains the support of the Fourier transform
$\widehat{f}$ of $f$. Put
$$
F(z)=\int_B \widehat{f}(u) e^{i(u, a)}e^{iz(u, b)}du,
\qquad z\in \mathbb{C}.
$$
Clearly $F(x)=f\circ l(x)$ for $x\in \mathbb {R}$. It remains to
note that since for $u\in B$ we have
$$
|e^{i(u, a)}e^{iz(u, b)}|=|e^{iz(u, b)}|\le e^{|z||(u, b)|}\le
e^{|z|r|b|},
$$
where $|b|$ is the length of $b$, it follows that
$$
|F(z)|\le e^{r|b||z|}\int_B |\widehat{f}(u)|du\le
e^{r|b||z|}|B|^{1/2}\|\widehat{f}\|_{L^2(\mathbb{R}^n)},
$$
where $|B|$ is the volume of $B$. The lemma is proved.

\quad

  Let us show that condition (i) of the theorem implies that
the map $\varphi$ is affine.

   For an arbitrary interval $I\subset\mathbb R$ let $1_I$ be
its characteristic function: $1_I(u)=1$ for $u\in I$, $1_I(u)=0$
for $u\not \in I$.

   Consider also the ``triangle'' function $\Delta_{(-2, 2)}$
supported on the interval $(-2, 2)\subset\mathbb R$, namely
$$
\Delta _{(-2, 2)}(u)=\max \bigg(1-\frac{|u|}{2} ,~0 \bigg),
\qquad u\in \mathbb R.
$$

  Denote by $^\vee$ the inverse Fourier transform.
By simple calculation we obtain
$$
(1_{(-1, 1)})^{^\vee}(x)=\frac{2\sin x}{x}, \quad
(1_{(0, 2)})^{^\vee}(x)=e^{ix}\frac{2\sin x}{x},
$$
$$
(\Delta_{(-2, 2)})^{^\vee}(x)=\frac{2\sin^2 x}{x^2}, \qquad x\in
\mathbb R,
$$
where we assume that
$$
\frac{\sin x}{x}\bigg|_{x=0}=1.
$$

   For $1\le j\le n$ consider the following functions
of $u=(u_1, u_2, \ldots, u_n)\in \mathbb{R}^n$ (products related
to empty set of indices are assumed to be equal to one)
$$
k(u)=\prod_{1\le s\le n} 1_{(-1, 1)}(u_s),
$$
$$
p_j(u)=\Delta_{(-2, 2)}(u_j)\prod_{\buildrel {1\le s\le n} \over
{s\neq j}} 1_{(-1, 1)}(u_s),
$$
$$
q_j(u)=1_{(0, 2)}(u_j)\prod_{\buildrel {1\le s\le n} \over
{s\neq j}} 1_{(-1, 1)}(u_s) .
$$
Define the functions $K, P_j, Q_j, ~j=1, 2, \ldots,n,$ on $\mathbb
R^n$ by
$$
K=(k)^{^\vee}, \quad P_j=(p_j)^{^\vee}, \quad Q_j=(q_j)^{^\vee}.
$$
We have
$$
K(t)=\prod_{1\le s\le n} \frac{2\sin t_s}{t_s},
$$
$$
P_j(t)= \frac{2\sin^2 t_j}{t_j^2}\prod_{\buildrel {1\le s\le n,}
\over {s\neq j}} \frac{2\sin t_s}{t_s},
$$
$$
Q_j(t)=e^{it_j}\frac{2\sin t_j}{t_j}\prod_{\buildrel {1\le s\le
n,} \over {s\neq j}} \frac{2\sin t_s}{t_s},
\qquad t=(t_1, t_2, \dots, t_n)\in\mathbb R^n.
$$

  We see that
$$
K, ~P_j, ~Q_j \in PW(\mathbb{R}^n). \eqno (7)
$$

  Note that
$$
Q_j(t)=e^{it_j}K(t), \qquad t\in\mathbb R^n, \eqno (8)
$$
and since
$$
P_j(t)=\frac{\sin
t_j}{t_j}K(t)=\frac{\frac{Q_j(t)}{K(t)}-\frac{K(t)}{Q_j(t)}}{2it_j}K(t)
$$
for $t_j, K(t), Q_j(t)\neq 0$, we have
$$
2iP_j(t)Q_j(t)t_j=Q_j(t)^2-K(t)^2, \qquad t\in\mathbb R^n. \eqno
(9)
$$

   We write the map $\varphi : \mathbb R^m\rightarrow\mathbb R^n$
in the form
$$
\varphi(t)=(\varphi_1(t), \varphi_2(t), \ldots, \varphi_n(t)),\qquad
t=(t_1, t_2, \ldots, t_m)\in\mathbb R^m.
$$
where the functions $\varphi_j :\mathbb R^m\rightarrow\mathbb R,
~j=1, 2, \ldots, n,$ are continuous.

   Let $l : \mathbb R\rightarrow \mathbb R^m$ be an arbitrary
affine map.

   Clearly the class $PW(\mathbb R^n)$ is invariant under
translations of the variable, so (see (7)) the functions
$$
K(t-\varphi\circ l(0)), ~P_j(t-\varphi\circ l(0)),
~Q_j(t-\varphi\circ l(0)), \qquad t\in\mathbb R^n,
$$
are in $PW(\mathbb R^n)$. From condition (i) of the theorem it
follows that the functions
$$
K(\varphi(t)-\varphi\circ l(0)), ~P_j(\varphi(t)-\varphi\circ l(0)),
~Q_j(\varphi(t)-\varphi\circ l(0)), \qquad t\in\mathbb R^m,
$$
are in $PW(\mathbb R^m)$.

  Consider the following functions on $\mathbb R$:
$$
K^*(x)=K(\varphi\circ l(x)-\varphi\circ l(0)),
$$
$$
P_j^*(x)=P_j(\varphi\circ l(x)-\varphi\circ l(0)),
$$
$$
Q_j^*(x)=Q_j(\varphi\circ l(x)-\varphi\circ l(0)), \qquad x\in\mathbb R.
$$
By Lemma 2 each of these functions is the restriction to $\mathbb
R$ of a corresponding entire function of exponential type.
Preserving notation we shall denote these functions by $K^*(z),
~P_j^*(z), ~Q_j^*(z)$.

  Note now that from (8) and (9) we have
$$
Q^\ast_j(x)=e^{i(\varphi_j\circ l(x)-\varphi_j\circ l(0))}
K^\ast(x), \qquad x\in\mathbb R,\eqno (10)
$$
and correspondingly
$$
2iP^\ast_j(x)Q^\ast_j(x)[\varphi_j\circ l(x)-\varphi_j\circ l(0)]=
(Q^\ast_j(x))^2-(K^\ast(x))^2,
\qquad x\in\mathbb R. \eqno (11)
$$

   Note also that
$$
K^\ast(0)=K(0)\not=0, \qquad
2iP^\ast_j(0)Q^\ast_j(0)=2iP_j(0)Q_j(0)\not=0,
$$
thus, the functions $K^\ast$ and $2iP^\ast_j Q^\ast_j$ are not
identically zero.

  We fix $j$. Clearly the sum and the product
of entire functions of exponential type are also entire functions
of exponential type. From  (10), (11), applying Lemma 1 to the
functions
$$
\psi=\varphi_j\circ l-\varphi_j\circ l(0), \quad
A=2iP^\ast_j Q^\ast_j, \quad
B=(Q^\ast_j)^2-(K^\ast)^2, \quad
C=K^\ast, \quad
D=Q^\ast_j,
$$
we obtain that $\varphi_j\circ l(x)-\varphi_j\circ l(0)$ is an
affine function of $x\in \mathbb{R}$. Therefore, $\varphi_j\circ
l(x)$ is an affine function of $x\in \mathbb{R}$. Since the affine
map $l : \mathbb R\rightarrow\mathbb R^m$ was chosen arbitrarily,
we obtain that the function $\varphi_j : \mathbb R^m \rightarrow
\mathbb R$ is affine. Since this is true for all $j, ~1\leq j\leq
n$, we see that $\varphi=(\varphi_1, \varphi_2, \ldots,
\varphi_n)$ is an affine map $\mathbb R^m\rightarrow\mathbb R^n$.

  We claim that the map $\varphi$ is injective.
Each function $F\in PW(\mathbb R^m)$ is the (inverse) Fourier
transform of a certain function in $L^1(\mathbb R^m)$, hence
$$
\lim_{|t|\rightarrow\infty} F(t)=0.
\eqno(12)
$$
Let $\varphi(t)=At+b$, where $A$ is an $n\times m$ matrix and
$b\in \mathbb{R}^n$. Suppose that $\ker A\neq \{0\}$. Take a
function $f\in PW(\mathbb R^n)$ with $f(b)\neq 0$. Let
$F=f\circ\varphi$. For all $t\in\ker A$ we have
$F(t)=f(At+b)=f(b)$. Thus, the superposition $F=f\circ\varphi$
does not satisfy (12). The implication (i)$\Rightarrow$(ii) is
proved.

   The proof of the implication (ii)$\Rightarrow$(i) is practically
trivial. Let $\varphi(t)=At+b$, where $A$ is an $n\times m$
matrix, $\ker A=\{0\}$, and $b\in\mathbb R^n$. We shall show that
for every function $f\in PW(\mathbb R^n)$ we have
$f\circ\varphi\in PW(\mathbb R^m)$.

  Since the map $\varphi$ is injective, we have
$n\geq m$. The simple case when $n=m$ has already been discussed
(see (1)). So, we can assume that $n>m$. Since the class $PW$ is
invariant under translations of the variable, we can also assume
that $b=0$, i.e., that $\varphi(t)=At$.

  Consider the image $\mathrm{im}\,\varphi$ of the map $\varphi$:
$$
\mathrm{im}\,\varphi=\{Ax : x\in \mathbb R^m\}.
$$
It is an $m$ -dimensional subspace of $\mathbb R^n$. Consider also
the following subspace $L$ of $\mathbb R^n$
$$
L=\{x=(x_1, x_2, \ldots, x_m, 0, \ldots, 0)\in\mathbb R^n\}.
$$
Let $S$ be the natural map of $\mathbb R^m$ onto $L$, namely
$$
S: (x_1, x_2, \ldots, x_m)\rightarrow
(x_1, x_2, \ldots, x_m, 0, \ldots, 0).
$$
Let $K$ be an invertible linear map of $\mathbb R^n$ onto itself
such that $K(\mathrm{im}\, \varphi)=L$. We put $Q=S^{-1}K\varphi$.
Then $Q$ is an invertible linear map of $\mathbb R^m$ onto itself.
We have $\varphi=K^{-1}SQ$ (by $S^{-1}, ~K^{-1}$ we denote the
maps inverse to $S$ and $K$ respectively).

   Thus, it suffices to verify that for every function
$f\in PW(\mathbb R^n)$ we have $f\circ S\in PW(\mathbb R^m)$. This
can be easily done as follows. Let $f\in PW(\mathbb R^n)$. Then
$$
f(Sx)=\int_{\mathbb R^n} \widehat{f}(u)e^{i(u, Sx)} du,
\qquad x\in\mathbb R^m.
$$
For $u=(u_1, u_2, \ldots, u_n)\in\mathbb R^n$ put
$$
u'=(u_1, u_2, \ldots, u_m), \quad u''=(u_{m+1}, u_{m+2}, \ldots u_n).
$$
We shall write $u=(u', u'')$. Consider a ball in $\mathbb R^n$
centered at $0$ that contains the support of the Fourier transform
$\widehat{f}$ of $f$. Let $r$ be the radius of this ball. We have
$$
f(Sx)=\int_{\mathbb R^m}
\bigg(\int_{\mathbb R^{n-m}}
\widehat{f}(u', u'')du''\bigg)e^{i(u', x)}du'=
\int_{\mathbb R^m}g(u')e^{i(u', x)}du',
$$
where
$$
g(u')=\int_{|u''|\leq \sqrt{r^2-|u'|^2}}\widehat{f}(u', u'')du''.
$$
The function $g$ vanishes outside of the ball centered at $0$ and
of radius $r$ in $\mathbb R^m$. At the same time (using the Cauchy
inequality) we obtain
$$
|g(u')|\leq \int_{|u''|\leq r}|\widehat{f}(u', u'')|du''\leq
c(r, n, m)\bigg(\int_{\mathbb R^{n-m}}|\widehat{f}(u', u'')|^2du''\bigg)^{1/2}
$$
(where $c(r, n, m)$ depends only on $r, n$ and $m$) and we see
that $g\in L^2(\mathbb R^m)$. The theorem is proved.

\quad

   I am grateful to S. V. Konyagin, who turned my attention
to the problem.

\begin{center}
\textsc{References}
\end{center}

\flushleft
\begin{enumerate}

\item  R. Paley, N. Wiener, \emph{Fourier transforms in the
    complex domain}, Amer. Math. Soc., New York, 1934.

\item  S. Azizi,  D. Cochran,  and J. N. McDonald, ``On the
    preservation of bandlimitedness under non-affine time
    warping'', Proc. of the 1999 Int. Workshop on Sampling
    Theory and Applications (SAMPTA), Aug. 11-14, 1999, Loen,
    Norway, The Norwegian University of Science and
    Technology, P. 37-40.

\item  S. Azizi,  J. N. McDonald, and  D. Cochran,
    ``Preservation of bandlimitedness under non-affine time
    warping for multi-dimensional functions'', In: 20th
    Century Harmonic Analysis -- A Celebration, J. S. Byrnes,
    ed., NATO Science Series, II Mathematics, Physics and
    Chemistry, 2001, V. 33, Kluwer, P. 369.

\item  E. C. Titchmarsh, \emph{The theory of functions}, 2nd
    edition, Oxford Univ. Press, New York, 1939.

\item  A. I. Markushevich, \emph{The theory of analytic
    functions}, in Russian, GITTL, ML., 1950.

\end{enumerate}

\quad

\qquad Dept. of Mathematical Analysis\\
\qquad Moscow State Institute of Electronics\\
\qquad and Mathematics (Technical University)\\
\qquad E-mail address: \emph{lebedevhome@gmail.com}

\end{document}